\documentclass{gtmon_a}
\pdfoutput=1

\usepackage[all]{xy}

\proceedingstitle{Proceedings of the School and Conference in Algebraic
Topology (The Vietnam National University, Hanoi, 9-20 August 2004)}
\conferencestart{9 August 2004}
\conferenceend{20 August 2004}
\conferencename{School and Conference in Algebraic Topology}
\conferencelocation{Vietnam National University, Hanoi, Vietnam}

\editor{John Hubbuck}
\givenname{John}
\surname{Hubbuck}

\editor{Nguy\~\ecircumflex{}n H V H\uhorn{}ng}
\givenname{H\uhorn{}ng}
\surname{Nguy\~\ecircumflex{}n}

\editor{Lionel Schwartz}
\givenname{Lionel}
\surname{Schwartz}

\title[Bocksteins and the nilpotent filtration on the cohomology of spaces]{Bocksteins and the nilpotent filtration\\ on the cohomology of spaces}

\author{Gerald Gaudens}
\givenname{Gerald}
\surname{Gaudens}
\address{Math\,Institut der Universit\"at Bonn\\
Beringstr. 1\\\newline
D-53115 Bonn\\Germany}
\email{gaudens@math.uni-bonn.de}

\urladdr{}

\volumenumber{11}
\issuenumber{}
\publicationyear{2007}
\papernumber{4}
\startpage{59}
\endpage{79}

\doi{}
\MR{}
\Zbl{}

\arxivreference{}

  
\received{12 January 2005}
\revised{12 April 2005}
\accepted{21 April 2005}
\published{14 November 2007}
\publishedonline{14 November 2007}
\proposed{}
\seconded{}
\corresponding{}
\editor{}
\version{}

\let\xysavmatrix\xymatrix
\def\xymatrix{\disablesubscriptcorrection\xysavmatrix}
\AtBeginDocument{\let\bar\wbar\let\tilde\wtilde\let\hat\what}
\AtBeginDocument{\def\R{\mathrm{R}}}
\def\cohom{\wwbar{\mathrm{H}}^{*}}
\def\redt{\wwbar{\mathrm{T}}}
\def\fin {\end{proof}}
\def \ie {ie~}

\makeop{ev}

\newdir{ >}{{}*!/-8pt/\dir{>}}
\newdir{ >}{{}*!/-5pt/@{>}}
\newdir{|>}{
    !/4.5pt/@{|}*:(1,-.3)@^{>}*:(1,+.3)@_{>}}

\def\tor{\mathrm{Tor}}

\def \lra {\longrightarrow}
\def \hra {\hookrightarrow}

\newcommand{\n}{\mathbb{N}}
\newcommand{\z}{\mathbb{Z}}

\makeatletter
\def\cnewtheorem#1[#2]#3{\newtheorem{#1}{#3}[section]
\expandafter\let\csname c@#1\endcsname\c@theoreme}

\newtheorem{theoreme}{Theorem}[section]
\cnewtheorem{lemme}[theoreme]{Lemma}
\cnewtheorem{proposition}[theoreme]{Proposition}
\cnewtheorem{cor}[theoreme]{Corollary}
\newtheorem*{conj}{Unbounded strong realization conjecture}

\theoremstyle{definition}
\cnewtheorem{definition}[theoreme]{Definition}
\cnewtheorem{remarque}[theoreme]{Remark}
\cnewtheorem{exm}[theoreme]{Example}
\cnewtheorem{exo}[theoreme]{Exercise}
\cnewtheorem{question}[theoreme]{Question}
\cnewtheorem{rem}[theoreme]{Remark}

\makeatother  %  move after \newtheorem block

\makeautorefname{lemme}{Lemma}
\makeautorefname{theoreme}{Theorem}

\def\lt{\mathrm{T}}

\begin{document}

\begin{webabstract}
N Kuhn has given several conjectures on the special features satisfied
by the singular cohomology of topological spaces with coefficients
in a finite prime field, as modules over the Steenrod algebra. The
so-called \emph{realization conjecture\/} was solved in special cases
in [Ann.  of Math. 141 (1995) 321--347] and in complete generality
by L Schwartz [Invent. Math. 134 (1998) 211--227]. The more general
\emph{strong realization conjecture\/} has been settled at the
prime $2$, as a consequence of the work of L Schwartz [Algebr. Geom.
Topol. 1 (2001) 519--548] and the subsequent work of F-X Dehon and
the author [Algebr. Geom. Topol. 3 (2003) 399--433]. We are here
interested in the even more general \emph{unbounded strong  realization
conjecture\/}. We prove that it holds at the prime $2$ for the class of
spaces whose cohomology has a trivial Bockstein action in high degrees.
\end{webabstract}

\begin{htmlabstract}
N Kuhn has given several conjectures on the special features satisfied
by the singular cohomology of topological spaces with coefficients
in a finite prime field, as modules over the Steenrod algebra. The
so-called <em>realization conjecture</em> was solved in special cases
in [Ann.  of Math. 141 (1995) 321&ndash;347] and in complete generality
by L Schwartz [Invent. Math. 134 (1998) 211&ndash;227]. The more general
<em>strong realization conjecture</em> has been settled at the
prime 2, as a consequence of the work of L Schwartz [<a
href="http://dx.doi.org/10.2140/agt.2001.1.519">Algebr. Geom.
Topol. 1 (2001) 519&ndash;548</a>] and the subsequent work of F-X Dehon and
the author [<a href="http://dx.doi.org/10.2140/agt.2003.3.399">Algebr.
Geom. Topol. 3 (2003) 399&ndash;433</a>]. We are here
interested in the even more general <em>unbounded strong  realization
conjecture</em>. We prove that it holds at the prime 2 for the class of
spaces whose cohomology has a trivial Bockstein action in high degrees.
\end{htmlabstract}

\begin{abstract}
N\,Kuhn has given several conjectures on the special features satisfied by
the singular cohomology of topological spaces with coefficients in a finite
prime field, as modules over the Steenrod algebra \cite{Ku}. The so-called
\emph{realization conjecture\/} was solved in special cases in \cite{Ku} and
in complete generality by L\,Schwartz \cite{Sc2}. The more general
\emph{strong realization conjecture\/} has been settled at the prime $2$, as a
consequence of the work of L\,Schwartz \cite{Sc3} and the subsequent work
of F-X\,Dehon and the author \cite{DG}. We are here interested in
the even more general \emph{unbounded strong  realization conjecture\/}. We
prove that it holds at the prime $2$ for the class of spaces whose
cohomology has a trivial Bockstein action in high degrees.
\end{abstract}

\maketitle

\section{Introduction}
The singular cohomology of a topological space with coefficients
in a finite prime field is naturally endowed with the structure of
an unstable algebra over the Steenrod algebra. That is, a graded ring
structure with a compatible action of the Steenrod algebra; see Schwartz
\cite[page 21]{Sc1}.

An unstable module isomorphic to the cohomology of some space is termed \emph{topologically realizable\/}.
N\,Kuhn's conjectures \cite{Ku} claim that realizable unstable modules
have rather special algebraic features. Namely, these conjectures tell us
that the action of the Steenrod algebra on the cohomology of a topological
space ought to be \emph{either very big or very small\/}.

The first of these conjectures \cite[Realization Conjecture, page 321]{Ku} was
settled by L\,Schwartz \cite[Theorem 0.1]{Sc2}  and says that the singular
cohomology of a space $X$ with coefficients in a finite prime field is
finitely generated as a module over the Steenrod algebra if and only if it
is finite dimensional as a (graded) vector space. In other words, the cohomology is nontrivial in finitely many degrees, and is a finite dimensional vector space in these degrees.

The more general \emph{strong realization conjecture\/} \cite[page 324]{Ku} was
settled at the prime $2$ by L\,Schwartz \cite{Sc3} under some
finiteness assumptions later removed by the work of F-X\,Dehon and the author in \cite{DG}.

Let us explain briefly the content of the \emph{strong realization conjecture\/}. Lannes' $\mathrm{T}$ functor is the endofunctor of the category of unstable modules which is left adjoint to tensoring with the cohomology of the infinite real projective space $\mathrm{H}^{*} \mathrm{B}(\z/2\z)$. There is a reduced version $\wwbar{\lt}$ of this 
functor which is left adjoint to tensoring with the \emph{reduced\/} cohomology of the infinite real projective space.
Let $\mathcal{U}_d$ be the full
subcategory of unstable modules  annihilated by $\redt ^{d+1}$,
the reduced Lannes' functor iterated $(d+1)$ times. The
subcategory $\mathcal{U}_ 0$ happens to be the subcategory of \emph{locally finite\/}
modules \cite{Sc1,Sc3}, \ie the full subcategory of unstable modules such that all monogenic submodules are finite dimensional over the ground field. The \emph{strong realization conjecture\/} says that if the singular
cohomology of a space $X$ with coefficients in a finite prime field is
in $\mathcal{U}_d$ for some $d$, then it is in $\mathcal{U}_0$.

From now on, we turn our attention to the even more general \emph{unbounded
strong  realization conjecture} (described in \fullref{theunb}), which we show to hold at the prime $2$ for the class of spaces having a trivial action of Bocksteins in high degrees.

\begin{rem}\label{rem11} In the following, $\mathrm{H}^* X$ always means the modulo $2$
singular cohomology of the space $X$. Also, all unstable modules are modules over the modulo $2$ Steenrod algebra.
\end{rem}

\subsection*{The main result}
\label{defnil}

We denote by $\mathcal{U}$ the category of unstable modules over the modulo $2$
Steenrod algebra. Every object $M$ of $\mathcal{U}$ is equipped with a natural
decreasing filtration, the so-called \emph{nilpotent filtration\/} \cite{Ku,Sc3}:
$$M= M_0\supset M_1 \supset  \ldots M_s \supset M_{s+1} \supset \ldots$$
This filtration is defined in the following way. An unstable module is called $s$--nilpotent if it belongs to the smallest full abelian subcategory of unstable modules containing $s$--th suspensions, and stable under extensions and filtered colimits.
The $s$--th step $M_s$ of the nilpotent filtration of an unstable module $M$ is its largest $s$--nilpotent submodule.

For each
$s$, the subquotients $M_s /M_{s+1}$ of the nilpotent filtration of $M$ are of the form $\Sigma ^s \mathrm{R}_s M$
where $\mathrm{R}_s M$ is a reduced module (see \fullref{nilpfilt}).

Our main result is the following.

\begin{theoreme}
\label{theo3}
Let $X$ be a topological space and let $\mathrm{H}^* X$ be its cohomology modulo $2$. Assume furthermore that $\mathrm{H}^* X$ has a trivial action of the Bockstein operator in high degrees and that $\mathrm{H}^* X$ is not locally finite.  The module $\mathrm{R}_t \cohom X$ cannot be locally finite for all integers $t\geq 0$, so let $s$ be the smallest $t$ such that $\mathrm{R}_t \cohom X$  is not locally finite. Then the unstable module $\mathrm{R}_s \cohom X$ does not belong to $\mathcal{U}_d$ for any integer $d$.
\end{theoreme}

The first assertion of the theorem follows from \fullref{uoparnil}.

So we get in particular that the \emph{unbounded strong realization conjecture\/} (to be explained in \fullref{theunb}) holds for the class of spaces such that the Bockstein acts trivially in high degrees:

\begin{theoreme}
\label{cor1}
Let $M$ be an unstable module such that for all $s$, the module $\mathrm{R}_s M$ is in some $\mathcal{U} _{d(s)}$. Suppose moreover that the Bockstein acts trivially on $M$ in high degrees. If $M$ is topologically realizable then $M$ is locally finite.
\end{theoreme}

In this statement, the number $d(s)$ is not supposed to be bounded with $s$; this explains the term \emph{unbounded\/} for the conjecture.
Let us explain briefly how \fullref{cor1} follows from \fullref{theo3}. Let $M$ be a topologically realizable unstable module $M$ such that the module $\mathrm{R}_s M$ is in some $\mathcal{U} _{d(s)}$ for all $s$ and such that $M$ has a trivial action of Bocksteins in high degrees. Suppose now, contradicting  \fullref{cor1}, that $M$ is not locally finite.
From \fullref{uoparnil}, we know that some $\mathrm{R}_s M$ is not locally finite. Assume $s$ is the smallest integer having this property. On the one hand, the hypotheses of \fullref{theo3} are fulfilled and $\mathrm{R}_s M$ is not 
in $\mathcal{U}_d$ for any $d$. But on the other hand, we had assumed the module $\mathrm{R}_s M$ to be in some $\mathcal{U} _{d(s)}$ for all $s$. This is a contradiction.

One might compare \fullref{cor1} to \cite[Theorem 0.1, Theorem 0.3]{Ku} in the seminal article of N\,Kuhn, where he proves the \emph{realization conjecture\/} under the
same hypothesis on Bocksteins as ours. The method he uses relies on secondary operations and does a priori
not apply to the more general setting of the \emph{unbounded strong realization conjecture\/}. We realized actually that the method of L\,Schwartz applies in our situation precisely in trying to extend (unsuccessfully) secondary operation technology to the more general realization conjectures.

Assume the unbounded conjecture is true in general (see \fullref{theunb}). If the cohomology ring $\cohom X$  of a space $X$ is not locally
constant, then for some integer $s$, the
reduced module $\mathrm{R}_s \cohom X$ does not belong $\mathcal{U}_d$ for any integer $d$.
L\,Schwartz has provided precise conjectures \cite[Conjecture 0.2, Conjecture
0.3]{Sc3} about the value of the smallest such $s$ in special cases. Our main theorem says
that in the case of the vanishing of Bocksteins in high degrees, the smallest $s$ such that $R_s\cohom X$ is not locally finite is also the smallest $s$ such that $R_s\cohom\!X$ does not belong to $\mathcal{U}_d$ for any integer $d$.

Example 0.11 of \cite[page 326]{Ku} is very useful in order to understand our result. Let $Y$ be the $s$--th bar filtration of
$\mathrm{B}\mathbb{C}\mathrm{P}^{\infty}$. Then $Y$  is a space with nilpotent cohomology (see \fullref{nilpfilt}) such that for $1\leq t <s$, the module $\mathrm{R}_t \cohom Y$ is in  $\mathcal{U}_t$, but $\mathrm{R}_s \cohom Y$ is not in $\mathcal{U}_d$ for any finite $d$. Our result shows that all cohomology classes which reduce nontrivially in $\mathrm{R}_1 \cohom  Y$ have a nonzero Bockstein.

This example shows that in general, if $\cohom X$ is not locally finite,
\begin{itemize}
\item the smallest value $s$ of $t$ such that  $\mathrm{R}_t \cohom X$ is not is not in $\mathcal{U}_d$ for any integer $d$ can be arbitrary high, 
\item the unstable modules $\mathrm{R}_t \cohom X$ for $1\leq t <s$ may be nonlocally finite. 
\end{itemize}

To prove \fullref{theo3}, we shall use, as in \cite{DG} the
theory of profinite spaces to be free of any finiteness hypotheses. The
\fullref{theo3} is a consequence of the more general:

\begin{theoreme}
\label{theo4}
Let $X$ be a \emph{profinite  space\/} and let $\mathrm{H}^* X$ be its continuous cohomology modulo $2$. Assume furthermore that $\mathrm{H}^* X$ has a trivial action of the Bockstein operator in high degrees and that $\mathrm{H}^* X$ is not locally finite.  The module $\mathrm{R}_t \cohom X$ cannot be locally finite for all integers $t\geq 0$, and we let $s$ be the smallest $t$ such that $\mathrm{R}_t \cohom X$  is not locally finite. Then the unstable module $\mathrm{R}_s \cohom X$ does not belong to $\mathcal{U}_d$ for any integer $d$.
\end{theoreme}

\fullref{theo4} implies \fullref{theo3} because the cohomology of a
space is naturally isomorphic to that of its profinite completion (which is
a profinite space) as an unstable algebra \cite[page 404, Section 2.3]{DG}.
Namely, suppose $X$ is a space such that $\mathrm{R}_s \cohom X$ is finite for
each $s$ and such that the Bockstein operator is zero in high degrees. Then the same holds for the cohomology of the profinite
completion of $X$. Hence, \fullref{theo4} implies \fullref{theo3}.

In the following, the word space means profinite space. Hence,
cohomology  means continuous cohomology, etc. What we need from the theory of profinite
spaces is strictly parallel to that of ordinary spaces.  All the constructions
on profinite spaces we will use are explained in detail in \cite{DG}. They behave in the same way as the usual constructions on spaces 
in the topological context. That's why the reader should not worry too much about profinite spaces in a first reading. From a philosophical viewpoint, profinite spaces are a replacement for usual spaces, where all our tools work without any restriction.

The setting of profinite spaces is crucial in the proofs,
for otherwise the tools we use (Lannes' functor, Eilenberg--Moore spectral sequence) would not work.

\section[Reformulations of the unbounded strong realization conjecture]{Reformulations of the  unbounded strong realization\\ conjecture}
\label{reform}

\subsection{Lannes' functor and the nilpotent filtration}
\label{nilpfilt}

The nilpotent filtration is briefly defined below \fullref{rem11}.

We begin by recalling an important property of the nilpotent filtration: any
unstable module is complete with respect to its nilpotent filtration. This
means that the natural map $M\lra \mathrm{lim}_s M/ M_s$ is an isomorphism.
This can be seen from the fact that for each $s$, the module $M_s$ is
$(s-1)$--connected.

We say that an unstable module is \emph{reduced\/} if the operator
$$\mathrm{Sq}_0\co M\lra M,~m\longmapsto \mathrm{Sq}^{|m|}m$$
is injective. If $M$ is the underlying module of some unstable algebra, then $M$
is reduced if and only this algebra has no nilpotent elements, because in any unstable algebra $M$, besides the \emph{Cartan formula\/} which says that
$$
\mathrm{Sq}^n(xy)= \sum _{i+j=n} (\mathrm{Sq}^i x) ~(\mathrm{Sq}^j y), \quad \text{for all } x,y \in M,
$$
we have the following other compatibly relation between the product and the Steenrod squares:
$$
\mathrm{Sq}_0 m = \mathrm{Sq}^{|m|}m = m^2 \quad \text{for all } m \in M. 
$$
In other words, the higher Steenrod square acting nontrivially, coincides with the Frobenius operator of the algebra.

For each
$s$, the subquotients $M_s /M_{s+1}$ of the nilpotent filtration of $M$ are of the form $\Sigma ^s \mathrm{R}_s M$
where $\mathrm{R}_s M$ is a reduced module.

On the other hand, an unstable module $M$ can be seen to be $1$--nilpotent (or simply nilpotent, for short) if and only if the operator $\mathrm{Sq}_0$ is locally nilpotent. This means that for all $m\in M$ there is a $t$ (depending a priori on $m$) such that
$$
(\mathrm{Sq}_0)^t m =0.
$$
An unstable module such that $M= M_s$ is called \emph{at least
$s$--nilpotent}. A $1$--nilpotent module is simply called \emph{nilpotent\/}. An
element of an unstable module is $s$--nilpotent provided it spans a
$s$--nilpotent submodule.

An important feature of the nilpotent filtration is its compatibility with tensor products: the tensor product of an $s$--nilpotent module with a $t$--nilpotent module is $(s+t)$--nilpotent.

The functor $\redt$ commutes with the nilpotent filtration in the
following sense (see \cite[Proposition 2.5, page 331]{Ku}):

\begin{proposition}
\label{comnilp}
Let $M$ be any unstable
module and let
$$
M= M_0\supset M_1 \supset \ldots M_s \supset M_{s+1} \supset \ldots \quad 
$$
be the nilpotent filtration of $M$. Then the induced filtration of $\redt M$
$$\redt M= \redt M_0\supset \redt M_1 \supset  \ldots \redt M_s
\supset \redt M_{s+1} \supset \ldots  \quad $$
is the nilpotent filtration of $\redt M$, \ie for all $s$,
$$\redt (M_s) =(\redt M)_s.$$
\end{proposition} 

As a consequence, by exactness and commutation of $\redt$ with suspensions, we have a sequence of equalities and natural isomorphisms \begin{align*}\Sigma ^s \mathrm{R}_s \redt M &=(\redt M)_s/(\redt M)_{s+1}\\&=\redt (M_s)/\redt (M_{s+1})\cong\redt (M_s /M_{s+1})=\redt \Sigma ^s \mathrm{R}_s
M \cong\Sigma ^s \redt \mathrm{R}_s M.\end{align*}
That is, the functors $\redt$ and $R_s$ commute for all $s$, up to natural
isomorphisms.

\subsection{Weight and the Krull filtration}
\label{weightdef}

Let $n$ be an integer. Let $n=\sum_{i=1}^{\ell} 2^{n_i}$ be the
binary expansion of $n$. We attach to $n$ the integer $\alpha (n)=\ell$.

\begin{definition}
\label{poids}
Let $M$ be a \emph{reduced\/} unstable module. We say that $M$ is of weight at most
$t$ if $M$ is trivial in all degrees $\ell$ such that $\alpha (\ell)>t$. The weight $w(M)$ of $M$ is the integer (maybe infinite) such that
$M$ is of weight at most $w(M)$ but not
$w(M)-1$.
\end{definition}

To understand the definition, we give the following examples.
\begin{exm}
Let $\mathrm{F}(1)$ be the unstable submodule generated by the nonzero
degree one class in $\mathrm{H}^{*} \mathrm{B}(\z/2\z)=\mathbb{F} _2 [u]$. It is
exactly the submodule of primitive elements of the Hopf algebra
$\mathrm{H}^{*} \mathrm{B}(\z/2\z)$. A graded $\mathbb{F}_2$--basis for
$\mathrm{F}(1)$ is given by the elements $\{u^{2^i}\}_{i\in\n}$. So
$\mathrm{F}(1)$ is zero in degrees $\ell$ such that $\alpha(\ell)$ is
strictly more that one. Hence the weight $w(\mathrm{F}(1))$ equals $1$.
\end{exm}
\begin{exm}
It is easy to see that $w(\mathrm{F}(1)^{\otimes n})= n$.
\end{exm}
\begin{exm}
The reduced cohomology ring $\cohom B(\z/2\z)=\mathbb{F} _2 [u]$ is of infinite
weight.
\end{exm}

A reduced module is of weight zero if and only if it is concentrated
in degree zero. In this case, we say that $M$ is \emph{constant\/}. For a reduced module, one readily checks that being \emph{constant\/} and \emph{locally finite\/} are equivalent notions.

More generally, the notion of weight and Krull filtration coincide for reduced modules, as shown by the following proposition.

\begin{proposition}{\rm (Franjou and Schwartz \cite{FS})}\qua
\label{comppoidsetkrull}
A reduced unstable module $M$ is in $\mathcal{U}_n$ if and only if its weight $w(M)$
is less or equal to $n$.
\end{proposition}

In particular, this  implies that a reduced module $M$ is in $\mathcal{U}_n$ if and only if $\redt ^n M \neq 0$ and $\redt ^{n+1} M= 0$. This proposition is an important tool for us, as we wish to consider the Krull filtration of the subquotients of the nilpotent filtration of certain unstable modules, and these subquotients are precisely reduced modules.

\subsection{The unbounded realization conjecture}
\label{theunb}

We can state the \emph{unbounded strong realization conjecture\/}
\cite[page 326]{Ku} in a slightly modified form.

\begin{conj}
\label{unbounded}
Let $M$ be an unstable module such that $\mathrm{R}_s M$ is of finite weight
for each $s$. If $M$ is topologically realizable, then the module
$\mathrm{R}_s M$ is constant for all $s$.
\end{conj}

The original conjecture of N\,Kuhn is not stated in terms of weight,
but in terms of polynomial degree of functors \cite[pages 325--326]{Ku}. This deserves a short explanation. Let $\mathcal{N}il$ be the full
subcategory of $\mathcal{u}$ of \emph{nilpotent\/} unstable modules. One can form the
quotient category $\mathcal{U}/\mathcal{N} il$. It is known  by Henn, Lannes and Schwartz \cite{HLS} that $\mathcal{U}/\mathcal{N}il$ is
equivalent to the full subcategory $\mathcal{F} _{\omega}$ of \emph{analytic
functors} of the category $\mathcal{F}$, where $\mathcal{F}$ is the category of functors
from finite dimensional
$\mathbb{F}_2$--vector spaces to all $\mathbb{F}_2$--vector spaces (with
natural transformations as morphisms). In the category $\mathcal{F}$, one has
a notion of polynomial functor of degree $n$.

Let $q\co\mathcal{U}\lra \mathcal{F} _{\omega} $ denote the quotient functor $
\mathcal{U}\lra\mathcal{U}/\mathcal{N} il$ composed with the equivalence of categories $\mathcal{U}/\mathcal{N} il \cong
\mathcal{F} _{\omega}$.

The point is that a reduced unstable module is of weight $n$ if and
only if $q(M)$ is polynomial of degree $n$.

We shall underline the proof of the
fact that the strong realization conjecture is a consequence of
the unbounded strong  realization conjecture. It relies on the following
lemma.
\begin{lemme}
\label{uoparnil}
An unstable module $M$ is in $\mathcal{U}_n$ if and only if $R_s M$ is in $\mathcal{U}_n$ for
all $s$.
\end{lemme}
\begin{proof}
Suppose $M$ is in $\mathcal{U}_n$. As $\mathcal{U}_n$ is a Serre subcategory (\ie abelian and stable under extensions \cite{Sc3}), the
modules $M_s$ and $M_s/M_{s+1}=\Sigma ^s R_s M$ are in $\mathcal{U}_n$ for each $s$. But the functor $\redt$ commutes
with suspensions and (more generally) with the nilpotent filtration (\fullref{comnilp}), so $R_s M$ is also in $\mathcal{U}_n$.

Conversely, if $R_s M$ is in $\mathcal{U}_n$ for all $s$, by exactness of
$\redt$ it follows that $M/M_s$ (recall that the nilpotent filtration is
decreasing) is in $\mathcal{U}_n$ for each $s$.
In other words,
$$\redt ^{n+1} (M)/ (\redt ^{n+1}(M))_s=\redt ^{n+1} (M)/ \redt
^{n+1}(M_s)\cong\redt ^{n+1} (M/M_s) =0 $$
for each $s$. But $\redt ^{n+1} (M)$ is complete with respect to its
nilpotent filtration, hence
$$\redt ^{n+1} (M)=0.$$
It follows that $M$ is in $\mathcal{U}_n$.
\fin

Now suppose we have an unstable module $M$ which is realizable and is
in $\mathcal{U}_n$, \ie such that $\redt ^{n+1} M =0$. By the preceding lemma, the
module $\mathrm{R}_s M$ is also in $\mathcal{U}_n$.
But an unstable module is of finite weight $n$ if and only if it is in
$\mathcal{U}_n$.

So, the unbounded strong  realization conjecture implies that
$\mathrm{R}_s M$ is constant for $s\geq 0$. Now, for a reduced module, being constant
and being in $\mathcal{U}_0$ are the same thing. Hence, by the lemma, the module $M$ is in
$\mathcal{u}_0$ and so the strong realization conjecture holds for $M$.

Another consequence of \fullref{uoparnil} is to give another form
of the unbounded strong realization conjecture:

\begin{conj}
\label{unbounded1}
Let $M$ be an unstable module such that $\mathrm{R}_s M$ is of finite weight
for each $s$. If $M$ is topologically realizable, then $M$ is locally
finite.
\end{conj}

This reformulation shows that \fullref{cor1} states a particular case of the \emph{unbounded strong realization conjecture\/}.

\section[Proof of Theorem 1.4]{Proof of \fullref{theo4}}

\subsection{Notations and summary of the proof}

It is not difficult to see that by replacing \emph{cohomology\/} by \emph{reduced cohomology\/} in \fullref{theo3}, one gets an equivalent statement. We will therefore work from now on with reduced cohomology.

The proof of \fullref{theo4} is by contradiction.
We want to prove that there exists no profinite space $X$ such that

\begin{enumerate} \item[(i)] the cohomology of $X$ is not
locally constant and for the lowest $d$ such that $\mathrm{R}_d
\cohom X$ is nonconstant, the module $\mathrm{R}_s \cohom X$ is of
finite weight, \item[(ii)] the action of the Bockstein is trivial
in high degrees in $\cohom X$.
\end{enumerate}
To this end, we refine the proof that was used in \cite{Sc1,Sc2,DG}. Let us recall how it goes.

Suppose that  a
profinite space $X$ satisfying the above conditions exists. Let $d$ be the
minimal integer $s$ such that $\mathrm{R}_s \cohom X$ is nonconstant.
Necessarily by \cite[Proposition 0.8, Corollary 0.9]{Ku}, $d$ is
nonzero. According to the discussion at the beginning of Section
7.2 in \cite{DG}, we can suppose that $\cohom X$ is $d$--nilpotent,
and as connected as necessary (the point here is that exchanging $X$ with the quotient of $X$ by some skeleton provides a new space with the same properties, but with higher connectivity).

We define for $0\leq\ell \leq d$, $$X_\ell=\Omega ^{d-\ell} X$$ so that
$X_d =X$ and $X_0 =\Omega ^d X$.

It follows from the
hypotheses that $\mathrm{R}_d \cohom X$ is of finite weight $f>0$. We
use \emph{Kuhn's reduction\/} in the framework of profinite spaces
\cite[Section 7.1]{DG} to lower the weight until $f=1$. This is done in \fullref{reduc}. This is the step that uses 
the technology of Lannes' $\mathrm{T}$ functor.

We construct a family ${(\alpha
_{i,d}})_{i\geq \kappa}$ of classes in $\cohom X$ satisfying a
certain set of conditions $(\mathcal{H}_d)$. We follow these
classes for $d\geq \ell \geq 0$ in the cohomology of the iterated
loop spaces $\Omega ^{d-\ell } X$: the classes
${(\alpha _{i,\ell }})_{i\geq \kappa}$ induced in $\Omega ^{d-\ell
} X$ through iterated evaluation map
$$\Sigma \Omega Z\lra Z$$
satisfy a
similar set of conditions $(\mathcal{H}_\ell )$. This is done in \fullref{weight} and \fullref{construc}. The properties of the Eilenberg--Moore spectral sequence are there heavily used.

The set of conditions $(\mathcal{H}_1)$ implies that the cup square of
$\alpha _{i,1}$ is trivial for large $i$ (see \fullref{cupnul1}). This is precisely the point where the hypothesis on the action of Bocksteins is needed. We show finally in \fullref{cupnul2}, following ideas of \cite{Sc3,DG} that the cup square of $\alpha _{i,0}$ is trivial for large $i$. Since the set of conditions $(\mathcal{H} _0)$ says in particular that the
cup square of $\alpha _{i,0}$ is nontrivial for large $i$, this
gives a contradiction.

An attentive reader may have noticed the method used here is very similar to that of \cite{Sc3,DG}. There are of course variations here, due the different situation. These are essentially
\begin{itemize}
\item we need to see that the hypothesis on Bocksteins carries over the Kuhn reduction (\fullref{reduc}),
\item the behaviour of the classes ${(\alpha _{i,d}})_{i\geq \kappa}$ is easier to analyse than in \cite{Sc3,DG}, because the set of hypotheses $(\mathcal{H}_\ell)$ is smaller,
\item we need on the other hand the slightly sharper statements on weight settled in \fullref{weight},
\item the last step explained in \fullref{cupnul2} is essentially the same as in \cite{Sc3,DG}, but in these sources, no clear statement we could rely on is made, and  the situation is also slightly different. We find it therefore useful to give full details in \fullref{cupnul2}.
\end{itemize}

\subsection{Kuhn's reduction with trivial Bocksteins}
\label{reduc}

Let $Y$ be a profinite space. Let $RY$ be the Bousfield--Kan functorial
fibrant replacement of $Y$ \cite{Mo} (see also \cite[Section 2.4]{DG}). We
denote by $\Delta Y$ the homotopy cofiber (in the homotopical algebra of
profinite spaces) of the natural map
$$Y\lra \mathrm{Map}(B(\z/2\z), RY).$$
Let $f\geq 1$ be the weight of $\mathrm{R}_d \cohom X$. We consider the
space $\Delta ^{f-1} X$.
\begin{lemme}
\label{redkuhn}
The space $\Delta ^{f-1} X$ satisfies
\begin{enumerate}
\item[(i)] the unstable module $\mathrm{R}_d \cohom \Delta ^{f-1} X$ is of
weight $1$,
\item[(ii)] the action of the Bockstein is trivial in high degrees in $\cohom \Delta ^{f-1} X$.
\end{enumerate}
\end{lemme}
\begin{proof}
It follows from \cite[Section 5]{DG} that
$$ \redt \cohom X \cong \cohom \Delta   X$$
as unstable modules.

As the nilpotent filtration commutes with $\redt$, it follows that
for all $s$ and $t$
$$\redt ^s \mathrm{R}_t \cohom X \cong \mathrm{R}_t \redt ^s\cohom X.$$
On the other hand, we know that $M$ is of weight $k$ if
and only if
$$ \redt ^{k+1} M=0 \quad \text{and} \quad \redt ^{k} M \neq 0.$$
We only need to prove that the action of the Bockstein is also
trivial in high degrees in $ \redt^{f-1} \cohom X \cong \cohom \Delta ^{f-1}   X$. But this is a
consequence of \fullref{appbock}.
\fin

\subsection{Weight watchers}
\label{weight}

We rely in this section and also in the last section on the properties of the Eilenberg--Moore spectral sequence for profinite spaces. We therefore recall the basic properties that will be used. Full details of its construction are given in \cite[Section 4]{DG}.

Let $X$ be a pointed profinite space. Then there is a natural second quadrant spectral sequence $\{(E_r^{-s, t}, d_r^{s,t}), s, t\geq0\}_{r\geq 1}$, converging to the cohomology of the loop space $\Omega X$, compatible with product and Steenrod operations. This means that for all $s\geq 0$ and $r\geq 2$, the graded vector space is an unstable module ${E}_{r}^{-s,*}$.
The differential
$$d_r\co {E}_{r}^{-s-r,*} \lra
{\Sigma}^{r-1} {E}_{r}^{-s,*}
$$
is linear with respect to the action of the Steenrod algebra.
The cohomology of the profinite loop space $\Omega X$ has a natural filtration by unstable submodules 
$$0= {\mathrm{F}}_0 \bar{\mathrm{H}}^* \Omega
X\subset {\mathrm{F}}_{-1} \bar{\mathrm{H}}^* \Omega X
\subset {\mathrm{F}}_{-2} \bar{\mathrm{H}}^* \Omega
X\subset \ldots \subset {\mathrm{F}}_{-s}
\bar{\mathrm{H}}^* \Omega X\subset \ldots
\subset \bar{\mathrm{H}}^* \Omega X$$
$${E}_{\infty}^{-s,*}\cong {\Sigma}^{s}({\mathrm{F}}_{-s}\bar{\mathrm{H}}^* \Omega X/{\mathrm{F}}_{-s+1}\bar{\mathrm{H}}^* \Omega X).
\leqno{\hbox{such that}}$$
This filtration converges to the cohomology of $\Omega X$ 
$${\bigcup}_{i\in\n} ~{\mathrm{F}}_{-i}\bar{\mathrm{H}}^* \Omega X =\bar{\mathrm{H}}^* \Omega X.$$
The spectral sequence carries products (in the most usual sense), and these products converge to the cup product on $\bar{\mathrm{H}}^* \Omega X$.

The $E_1$--term is given by the bar construction (see Mac Lane \cite{Mac}) and in particular  $E_1^{-s, t}= (\mathrm{H}^* X)^{\otimes s}$. The product on the $E_1$--term is given by the shuffle product \cite{Mac} and the Steenrod module structure is the canonical one. Thus the $E_2$--term is given by
$$E_2^{-s, t} = \tor^{-s, t}_{H^* X}(\mathbb{F}_ 2,\mathbb{F}_ 2).$$
No finiteness hypotheses are needed here to analyse the $E_2$--term as a $\tor$ group because we use the profinite setting \cite{DG}.

With the help of the Eilenberg--Moore spectral sequence, we will prove the following lemma.
\begin{lemme}
\label{R1} For $1\leq \ell \leq d$, the module $\mathrm{R}_\ell
\cohom X_\ell$ has weight one.
\end{lemme}
\begin{proof}
If $d=1$ the lemma is clearly true from the hypotheses, otherwise we prove
\fullref{R1} by induction on:
\begin{lemme}
\label{R12}
Let $Y$ be a profinite space such that $\cohom Y$ is $h$--nilpotent, $h\geq 2$.
Then $\mathrm{R}_{h-1} \cohom \Omega Y$ and $\mathrm{R}_h Y$ have the same
weight.
\end{lemme}
\proof[Proof of \fullref{R12}]
We use the Eilenberg--Moore spectral sequence which calculates $\cohom \Omega Y$
from $\cohom Y$. Its $E_{2}^{-s,*}$--term is a subquotient of $\smash{{\cohom Y}^{\otimes
s}}$, which is $hs$--nilpotent (because of the compatibility of tensor products with nilpotency, see \fullref{nilpfilt}). Because the subcategory of $t$--nilpotent
modules is a Serre subcategory (\ie abelian and stable under extensions), it happens that $E_{\infty}^{-s,*}$ is also $sh$--nilpotent.

Let $\{\mathrm{F}_{-s} \cohom \Omega Y\}_{s\in \mathbb{N}}$ be the Eilenberg--Moore filtration,
whose associated graded is the abutment of the Eilenberg--Moore spectral
sequence. We have
$$E_{\infty}^{-s,*}=\Sigma ^s (\mathrm{F}_{-s}/\mathrm{F}_{-s+1}) \cohom \Omega
Y$$
as unstable modules, hence $(\mathrm{F}_s/\mathrm{F}_{s-1}) \cohom \Omega Y$ is
$(hs-s)$--nilpotent.

Because the Eilenberg--Moore filtration is convergent, and $s$--nilpotent modules form a Serre subcategory stable under filtered colimits, we have that $\cohom \Omega Y/\mathrm{F}_s \cohom \Omega Y $ is at least $(hs-s)$--nilpotent.

We recall  the following result \cite[Corollary A.3]{DG}.

\begin{proposition}
\label{U-Nil}
Let $0\to A\to B\to C\to 0$ be a short exact sequence of unstable modules and $p,q,s$ three nonnegative integers. Suppose that $\R_sA$ is in $\mathcal{U}_p$ and that $\R_sC$ is in $\mathcal{U}_q$ ; then $\R_sB$ is in $\mathcal{U}_{\mathrm{max}\{p,q\}}$.
\end{proposition}

Applying   this result to the short exact sequence
$$
\mathrm{F}_{-1}\cohom \Omega Y \lra \cohom \Omega Y \lra \cohom \Omega
Y/\mathrm{F}_{-1} \cohom \Omega Y 
$$
we easily get that $\mathrm{R}_{h-1} \cohom \Omega Y$ and $\mathrm{R}_{h-1} \mathrm{F}_{-1} \cohom \Omega Y$ have the same weight.

We know that
$$
\mathrm{R}_{h-1} \mathrm{F}_{-1} \cohom \Omega Y\cong
\mathrm{R}_{h} \Sigma (\mathrm{F}_{-1} \cohom \Omega Y/\mathrm{F}_{0} \cohom \Omega Y)=\mathrm{R}_{h}E_{\infty}^{-1,*}
$$
and so we need to compare $\mathrm{R}_{h}E_{\infty}^{-1,*}$ and
$\mathrm{R}_{h} \cohom Y$.

But $E_{\infty}^{-1,*}$ is isomorphic to the quotient of
$\cohom Y$ by $B$, the union of the images of all higher differentials. The image of
the differential $d^r$ is easily seen to be at least
$((r+1)(h-1)+2)$--nilpotent, by using the linearity of differentials (see \cite{DG} for more details). Hence, the union of the image of
the differentials is at least $(2h-1)$--nilpotent. We have a short exact
sequence:
$$ B\lra \cohom Y \lra E_{\infty}^{-1,*}$$
A new application of \fullref{U-Nil} gives that $\mathrm{R}_h
E_{\infty}^{-1,*}$ and $\mathrm{R}_h \cohom Y$ are of the same weight, and
\fullref{R12} follows.
\fin
\begin{lemme}
\label{lemmefinal}
The module $\mathrm{R}_0 \mathrm{F}_{-1} \cohom X_0$ is of weight $1$. The
module $\mathrm{R}_0 \mathrm{F}_{-2} \cohom X_0$ is of weight $2$.
\end{lemme}

\begin{proof}
We have isomorphisms
$$\mathrm{R}_0 \mathrm{F}_{-1} \cohom X_0 \cong \mathrm{R}_1 \Sigma
(\mathrm{F}_{-1}\cohom X_0)/\mathrm{F}_{0} \cohom X_0)\cong \mathrm{R}_1
E_{\infty}^{-1, *}.$$
The module $E_{\infty}^{-1, *}$ is a quotient of $\cohom X_1$ by an at
least $2$--nilpotent submodule $B$.

So we have an exact sequence
$$B\lra \cohom X_1 \lra E_{\infty}^{-1, *}.$$
By \fullref{R1}, the module $\mathrm{R}_1 \cohom X_1$ is of weight
$1$ which proves the first assertion.

The module $\mathrm{R}_0 (\mathrm{F}_{-2}/\mathrm{F}_{-1}) \cohom X_0$
is isomorphic to $\mathrm{R}_2 \Sigma^2 (\mathrm{F}_{-2}/\mathrm{F}_{-1})
\cohom X_0 = \mathrm{R}_2 E_{\infty}^{-2, *}$. The module $E_{\infty}^{-2, *}$
is a subquotient of $(\cohom X_1 )^{\otimes 2}$. So we have modules $B\subset C
\subset (\cohom X_1 )^{\otimes 2}$ such that $C/B
=E_{\infty}^{-2, *}$. The module $B$ is the union of all the images of the
differentials and $C$ is the submodule of infinite cycles. One estimates
that $B$ is at least $3$--nilpotent.
Hence by \cite[Corollary A.2]{DG} implies that $\mathrm{R}_2
E_{\infty}^{-2, *}$ is isomorphic to $\mathrm{R}_2 C $. On the other hand
the functor $\mathrm{R}_2$ preserves monomorphisms \cite[Proposition A.1]{Sc3,DG} and so $\mathrm{R}_2 E_{\infty}^{-2, *}$ is isomorphic to some
submodule of $\mathrm{R}_2 ((\cohom X_1 )^{\otimes 2}$. We finally note that
$$\mathrm{R}_2 (\cohom X_1 )^{\otimes 2} = {\oplus}_{i+j=2} \mathrm{R}_i (\cohom X_1 )\otimes   \mathrm{R}_j (\cohom X_1 )=\mathrm{R}_1 (\cohom X_1 )\otimes
\mathrm{R}_1 (\cohom X_1 ).$$
As $\mathrm{R}_1 (\cohom X_1 )$ is of weight one, the module $\mathrm{R}_2 (\cohom X_1 )^{\otimes 2}$ is of weight $2$, and so are $\mathrm{R}_2
E_{\infty}^{-2, *}$ and $\mathrm{R}_0 \mathrm{F}_{-1} \cohom X_0$.
Using the short exact sequence
$$\mathrm{F}_{-1} \cohom X_0 \lra \mathrm{F}_{-2} \lra (\mathrm{F}_{-2} \cohom X_0/\mathrm{F}_{-1})\cong \Sigma^2 E_{\infty}^{-2, *}$$
and applying \fullref{U-Nil} and the preceding remarks, we find that
the module $\mathrm{R}_0 \mathrm{F}_{-2} \cohom X_0$ is of weight $2$.
\fin
\subsection{Construction of classes}
\label{construc}

The next lemma is a special case of Proposition 7.2 of Dehon and the author \cite{DG}. The original statements are in Schwartz \cite{Sc2,Sc3}.
\begin{lemme}
\label{u1} Let $M$ be a reduced module of weight $1$. Let $\eta$
be the unity of the adjunction $M\to \wwbar{\lt}M\otimes
\bar{\mathrm{H}}^{*}\mathrm{B}(\z/2\z)$. Then $\eta$ factorizes by
the submodule $\wwbar{\lt}  M \otimes \mathrm{F}(1)$. Moreover,
the kernel and cokernel of $$\eta\co M\to \wwbar{\lt}M\otimes
\mathrm{F}(1)$$ are locally finite.
\end{lemme}
We apply this lemma to $M=\mathrm{R}_d \cohom X$, which we can suppose to be
of weight $1$ by \fullref{redkuhn}. Then it follows that
there is a cyclic submodule of the form $\smash{\mathrm{F}(1)^{\geq 2^\xi}}$
in $M$, generated by some $\bar{\alpha} _\xi$ of degree $2^\xi$.
We can suppose $\xi$ as big as we want. So we pick up some
$\kappa\geq \xi$.

We lift up $\Sigma ^s \bar{\alpha}_\kappa$ to a class $\alpha
_{\kappa, d}$ of degree $2^\kappa +d$ through the epimorphism $(\cohom X)_s
\lra \Sigma ^s \mathrm{R}_s (\cohom X)$, and we define recursively, for $i\geq
\kappa$
$$ \alpha _{i+1,d}=\mathrm{Sq}^{2^i}\alpha _{i,d}.$$
We get some classes $(\alpha _{i,d})_{i\geq \kappa}$ satisfying the
following set of conditions:
$$
({\mathcal{H}}_d)\left\{\begin{array}{l}
\mathrm{the~class~}\alpha _{i,d}\mathrm{~is~defined~for~}i\geq
\kappa\mathrm{~and~is~of~degree~}2^ i +d\mathrm{~in~} \cohom X,\\
\mathrm{the~ class~ }\alpha _{i,d}\mathrm{ ~reduces ~ nontrivially ~in~}
\mathrm{R}_d (\cohom X)~\mathrm{ ~(hence ~ is ~ nonzero)},\\
\mathrm{the ~Bockstein ~acts ~trivially ~on~ }\alpha _{i,d},\\
\mathrm{for~}i\geq \kappa\mathrm{,~ we~ have~}\smash{\mathrm{Sq}^{2^i}}\alpha
_{i,d}=\alpha_{i+1,d}.
\end{array}\right.
$$
The evaluation $\ev _Z \co \Sigma \Omega Z \lra Z$ induces a morphism
$$
\ev _Z \co \mathrm{H}^* Z \lra\mathrm{H}^* \Sigma \Omega Z \cong \Sigma \mathrm{H}^* \Omega Z.
$$
We define iteratively, for all $0\leq \ell \leq d-1$, the classes $(\alpha _{i,\ell})_{i\geq \kappa}$ as
$$
\alpha _{i,\ell} = \ev _{X_{l+1}} \alpha _{i,\ell +1}.
$$
We prove by downward induction the following proposition.

\begin{proposition}
\label{propclass}
The classes $(\alpha _{i,\ell})_{i\geq \kappa}$ satisfy, for $0\leq\ell \leq
d$ and $i\geq \kappa$:
$$({\mathcal{H}}_\ell) \left\{\begin{array}{l}
\mathrm{the~class~}\alpha _{i,\ell}\mathrm{~is~of~degree~}2^ i
+\ell\mathrm{~in~} \cohom X,\\
\mathrm{the~ class~ }\alpha _{i,\ell}\mathrm{ ~reduces ~ nontrivially ~in~}
\mathrm{R}_\ell (\cohom X),\\
\mathrm{the ~Bockstein ~acts ~trivially ~on~ }\alpha _{i,\ell},\\
\mathrm{for~}i\geq \kappa,\mathrm{~ we~ have~}\smash{\mathrm{Sq}^{2^i}}\alpha
_{i,\ell}=\alpha_{i+1,\ell}.
\end{array}\right.$$
\end{proposition}
\begin{proof} The assertion on the degree of
$(\alpha _{i,\ell})_{i\geq \kappa}$ follows from the definitions. The second
point is a consequence of the following lemma (see \cite[Proposition
A.4]{DG}).
\begin{lemme}
Let $Y$ be a profinite space such that $\cohom Y$ is $\ell$--nilpotent for $\ell
\geq 1$. Then $\cohom Y$ is $(\ell -1)$--nilpotent and the evaluation morphism
induces a \emph{monomorphism\/}
$$\mathrm{R}_d\cohom Y  \hra \mathrm{R}_d\Sigma \cohom \Omega Y\cong
\mathrm{R}_{d-1} \Omega Y.$$
\end{lemme}
The third and fourth points are consequences of the Steenrod algebra
linearity of the evaluation morphism. Namely, it follows from the equalities
\begin{eqnarray*}
\Sigma (\mathrm{Sq}^1 \alpha _{i,\ell-1})&=&\mathrm{Sq}^1 \Sigma  \alpha _{i,\ell-1}\\
&=&
\mathrm{Sq}^1 \ev_{X_\ell}~ (\alpha _{i,\ell})\\
&=&
\ev_{X_\ell}~ (\mathrm{Sq}^1 \alpha _{i,\ell})\\
&=&
0
\end{eqnarray*}
that the Bockstein acts trivially on $\alpha _{i,\ell}$, and the equalities
\begin{eqnarray*}
\Sigma (\mathrm{Sq}^{2^i} \alpha _{i,\ell-1})
&=&
\mathrm{Sq}^{2^i} \Sigma \alpha _{i,\ell}\\
&=&
\mathrm{Sq}^{2^i} \ev_{X_\ell} (\alpha _{i,\ell})\\
&=&
\ev_{X_\ell} (\mathrm{Sq}^{2^i}\alpha _{i,\ell})\\
&=&
\ev_{X_\ell}~(\alpha _{i+1,\ell})\\
&=&
\Sigma \alpha _{i+1,\ell-1}
\end{eqnarray*}
show how $\mathrm{Sq}^{2^i}$ acts on $\alpha _{i,\ell}$. \fin

\subsection[The cup square of alpha sub (i,1) is trivial]{The cup square of $\alpha _{i,1}$ is trivial}
\label{cupnul1}

This is exactly the point where the hypothesis that Bocksteins are trivial in high degrees is used.

For $\ell=1$, the classes $\alpha _{i,1}$
have degree $2^ i +1$, and the unstable algebra structure gives for $i\geq
\kappa$,
$$\alpha _{i,1}\cup \alpha _{i,1}= \mathrm{Sq}^{2^ i +1}\alpha _{i,1}=
\mathrm{Sq}^1 \mathrm{Sq}^{2^ i}\alpha _{i,1}= \mathrm{Sq}^1 \alpha
_{i+1,1}=0.$$
So to sum up the situation, we have a profinite space $X_{1}=\Omega
^{d-1} X$ and classes
${(\alpha _{i,1})}_{i\geq \kappa}$ such that for $i\geq \kappa$,

\begin{enumerate}
\item[(i)] the class $\alpha _{i,1}$ is of degree $2^ i +1$ in $\cohom X_1$,
\item[(ii)] the class $\alpha _{i,1}$ reduces nontrivially in $\mathrm{R}_1
(\cohom X_1)$,
\item[(iii)] the Bockstein acts trivially on $\alpha _{i,1}$,
\item[(iv)] we have $\smash{\mathrm{Sq}^{2^i}}\alpha _{i,1}=\alpha _{i+1,1}$,
\item[(v)] the cup square $\alpha _{i,1}\cup \alpha _{i,1}$ is trivial.
\end{enumerate}

Suppose that we are able to prove that the same set of
conditions holds for ${(\alpha _{i,0})}_{i\geq \kappa'}$, then we
obtain the following contradiction
$$0=\alpha _{i,0}\cup \alpha _{i,0}=\mathrm{Sq}^{2^i}\alpha _{i,0}=\alpha
_{i+1,0}\neq 0.$$
So we need to prove that $\alpha _{i,0}\cup \alpha _{i,0}=0$,
for $i\geq \kappa$.

\subsection[The cup square of alpha sub (i,0) is trivial]{The cup square of $\alpha _{i,0}$ is trivial}
\label{cupnul2}

We use the Eilenberg--Moore spectral sequence which relates $\cohom X_{1}$ to $\cohom X_0=\cohom \Omega X_1$. Recall that the Eilenberg--Moore spectral sequence carries products in the following way: the shuffle product $[\,\cdot\,|\,\cdot\,]$ on the $E_1$--term of the Eilenberg--Moore spectral
sequence converges to the cup product on the $E_{\infty}$--term (which means in particular that the shuffle product of infinite cycles is itself an infinite cycle).

For $i\geq
\kappa$, the cup square $\alpha _{i,1}\cup \alpha _{i,1}$ is
trivial. So the element $\alpha _{i,1}\otimes \alpha _{i,1}= [\alpha _{i,1}, \alpha _{i,1}]$ is a $1$--cycle and defines an element of $E_{2}^{-1, *}$, as $\alpha _{i,1}\cup \alpha _{i,1}= d_1 (\alpha _{i,1}\cup \alpha _{i,1})$. For degree reasons, the higher differentials
coming from $\smash{E_{2}^{-1, *}}$ are trivial and so, the $1$--cycle
$\alpha _{i,1} \otimes \alpha _{i,1}$ induces a permanent cycle, which
never bounds for nilpotence reasons (see \cite[Section 7.4]{DG}).
Let $w_{i,\ell}$ be any element of $\cohom X_0$ detected by this
permanent cycle.

{\bf First step}\qua We want to compare
$\mathrm{Sq}^{2^i}w_{i,0}$ to $\alpha _{i+1,0}\cup \alpha
_{i,0}$. The cycle $[\alpha_{i,1}\,|\,\alpha _{i+1,1}]=\alpha_{i,1}\otimes \alpha _{i+1,1}+\alpha _{i+1,1}\otimes\alpha _{i,1}$ detects the cup product $\alpha _{i+1,0}\cup\alpha _{i,0}$.

By Cartan's formula, we have
$$\mathrm{Sq}^{2^i}[\alpha_{i,1}\,|\,\alpha _{i,1}]=\mathrm{Sq}^{2^i}(\alpha _{i,1}\otimes \alpha _{i,1})=
[\alpha _{i,1}|\,\alpha _{i+1,1}]+
{\sum}_{0<t\leq 2^{i-1}}[{\mathrm{Sq}}^{t}\alpha _{i,1}\,|\,{\mathrm{Sq}}^{2^i-t}\alpha _{i,1}].$$
The permanent cycle $\mathrm{Sq}^{2^i}(\alpha _{i,1}\otimes \alpha
_{i,1})$ converges to $\mathrm{Sq}^{2^i}w_{i,0}$ by compatibility
of the Eilenberg--Moore spectral sequence with Steenrod operations. In the same way, $[\mathrm{Sq}^{t}\alpha _{i,1}\,|\,\smash{\mathrm{Sq}^{2^i-t}}\alpha _{i,1}]$ converges to ${\mathrm{Sq}}^{t}\alpha _{i,1}\cup \smash{{\mathrm{Sq}}^{2^i-t}}\alpha _{i,1}$ for $0\leq t\leq 2^{i-1}$.

Therefore, the element 
$$
\mathrm{Sq}^{2^i}w_{i,0}-\alpha _{i,0}\cup
\alpha _{i+1,0}-{\sum}_{0<t\leq 2^{i-1}}{\mathrm{Sq}}^{t}\alpha _{i,0}\cup{\mathrm{Sq}}^{2^i-t}\alpha _{i,0}
$$
is in $\mathrm{F}_{-1}\cohom X_d$. This equation is homogeneous of degree $2^i+2^{i+1}$ and
$\alpha(2^i+2^{i+1})=2$ (the function $\alpha$ is defined in the beginning of \fullref{weightdef}). But $\mathrm{R}_0 \mathrm{F}_{-1}\cohom X_d$
is of weight $1$, by \fullref{lemmefinal}, so 
$$
\mathrm{Sq}^{2^i}w_{i,0}-\alpha _{i,0}\cup
\alpha _{i+1,0}-{\sum}_{0<t\leq 2^{i-1}}{\mathrm{Sq}}^{t}\alpha _{i,0}\cup{\mathrm{Sq}}^{2^i-t}\alpha _{i,0}
$$
reduces to zero in $\mathrm{R}_0\mathrm{F}_{-1}\cohom X_d$. Moreover,
the inclusion of  $\mathrm{F}_{-1}\cohom X_d$ in $\mathrm{F}_{-2}\cohom X_d$ induces a monomorphism $\mathrm{R}_0\mathrm{F}_{-1}\cohom X_d\lra \mathrm{R}_0\mathrm{F}_{-2}\cohom X_d$ \cite[Proposition
A.1]{DG}. Hence $$\mathrm{Sq}^{2^i}w_{i,0} \quad\text{and}\quad \alpha
_{i+1,0}\cup \alpha _{i,0}+{\sum}_{0<t\leq 2^{i-1}}{\mathrm{Sq}}^{t}\alpha _{i,0}\cup{\mathrm{Sq}}^{2^i-t}\alpha _{i,0}$$ project to \emph{equal\/} elements in $\mathrm{R}_0\mathrm{F}_{-2}\cohom X_d$.

We now note, that for $0<t<2^{i-1}$, the class $\mathrm{Sq}^{t}\alpha _{i,0}$ is in degree $t+2^i$ and $\alpha(t+2^{i})=2$. Therefore ${\mathrm{Sq}}^{t}\alpha _{i,0}$ reduces to zero in $\mathrm{R}_0 \mathrm{F}_{-1}\cohom X_d$, which is of weight $1$ by \fullref{lemmefinal}. Now, the product map
$$
{F}_{-1}\cohom X_d \otimes {F}_{-1}\cohom X_d \lra {F}_{-2}\cohom X_d
$$
induces a map
$$
(\mathrm{R}_0 {F}_{-1}\cohom X_d )\otimes (\mathrm{R}_0 {F}_{-1}\cohom X_d)\cong\mathrm{R}_0 ({F}_{-1}\cohom X_d \otimes {F}_{-1}\cohom X_d) \lra \mathrm{R}_0 {F}_{-2}\cohom X_d.
$$
It follows that $\textstyle{\sum}_{0<t\leq 2^{i-1}}{\mathrm{Sq}}^{t}\alpha _{i,0}\cup{\mathrm{Sq}}^{2^i-t}\alpha _{i,0}$ reduces to zero in $\mathrm{R}_0 {F}_{-2}\cohom X_d $, hence the following lemma holds:

\begin{lemme}
The elements $\mathrm{Sq}^{2^i}w_{i,0}$ and $\alpha _{i+1,0}\cup \alpha _{i,0}$ project to \emph{equal\/} elements in $\mathrm{R}_0\mathrm{F}_{-2}\cohom X_d$.
\end{lemme}

The class $\alpha_{i+1,0}\cup \alpha _{i,0}$ is in degree $2^i+2^{i+1}$ and $\alpha(2^i+2^{i+1})=2$.

{\bf Second step}\qua We now proceed to compare
$\mathrm{Sq}^{2^i}\mathrm{Sq}^{2^i}w_{i,0}$ and $\mathrm{Sq}^{2^i}(\alpha
_{i+1,0}\cup \alpha _{i,0})$.

The Cartan formula gives
\begin{align*}
\mathrm{Sq}^{2^i}( \alpha _{i,0}\cup \alpha _{i+1,0})&=
{\sum}_{p+q=2^i}\mathrm{Sq}^{p}\alpha _{i,0}\cup \mathrm{Sq}^{q}\alpha
_{i+1,0}\\
\mathrm{Sq}^{2^i} (\alpha _{i,0}\cup \alpha _{i+1,0})&=
(\mathrm{Sq}^{2^i}\alpha _{i,0})\cup \alpha _{i+1,0}+
{\sum}_{p<2^i}\mathrm{Sq}^{p}\alpha _{i,0}\cup \mathrm{Sq}^{2^i -p}\alpha
_{i+1,0}\tag*{\hbox{so that}}\\
&=\alpha _{i+1,0}\cup \alpha _{i+1,0}+ {\sum}_{p<2^i}\mathrm{Sq}^{p}\alpha
_{i,0}\cup \mathrm{Sq}^{2^i -p}\alpha _{i+1,0}.
\end{align*}
For $p<2^i$, we have the following two cases.
\begin{itemize}
\item If $0<p<2^i$, then $\mathrm{Sq}^{p}\alpha _{i,0}$ has degree
$\ell= 2^i +p$ such that $\alpha (\ell)>1$. The element $\alpha _{i,0}$ is
in the submodule $\smash{\mathrm{F}_{-1}\cohom X_0}$ by definition, thus so is
$\mathrm{Sq}^{p}\alpha _{i,0}$. But $\mathrm{R}_0
\mathrm{F}_{-1}\cohom X_0$ is of weight one and this implies that
$\mathrm{Sq}^{p}\alpha _{i,0}$ reduces to zero in $\mathrm{R}_0
\mathrm{F}_{-1}\cohom X_0$. In other words, the class
$\mathrm{Sq}^{p}\alpha _{i,0}$ is nilpotent.
\item If $p=0$, then $\mathrm{Sq}^{2^i -p}\alpha _{i+1,0}$ has degree
$\ell$ such that $\alpha (\ell)>1$. The same argument shows that if $p=0$,
the element $\smash{\mathrm{Sq}^{2^i -p}}\alpha _{i+1,0}$ is nilpotent.
\end{itemize}
So for $p<2^i$, either $\mathrm{Sq}^{p}\alpha _{i,0}$ or
$\smash{\mathrm{Sq}^{2^i -p}}\alpha _{i+1,0}$ is nilpotent and so is the cup product
$\mathrm{Sq}^{p}\alpha _{i,0}\cup \smash{\mathrm{Sq}^{2^i -p}}\alpha _{i+1,0}$.

It follows that $\smash{\mathrm{Sq}^{2^i}}( \alpha _{i,0}\cup \alpha
_{i+1,0})$ and $\alpha _{i+1,0}\cup \alpha _{i+1,0}$ project to \emph{equal\/}
elements in $\mathrm{R}_0\mathrm{F}_{-2}\cohom X_d$. In other words the following lemma holds.

\begin{lemme}
The classes $\mathrm{Sq}^{2^i}\mathrm{Sq}^{2^i}w_{i,0}$ and $\alpha _{i+1,0}\cup \alpha
_{i+1,0}$ project to equal elements in $\mathrm{R}_0\mathrm{F}_{-2}\cohom X_d$.
\end{lemme}

The decomposition of
$\mathrm{Sq}^{2^i}\mathrm{Sq}^{2^i}$ \cite[Lemma 5.7, page 554]{Sc3} implies
that  $\smash{\mathrm{Sq}^{2^i}\mathrm{Sq}^{2^i}}w_{i,0}$  belongs to a submodule
of $\mathrm{F}_{-2}\cohom X_0$ generated by elements having degrees $\ell$ such
that $\alpha(\ell)\geq 3$. But $\mathrm{R}_0 \mathrm{F}_{-2}\cohom X_0$ is of
weight $2$ by \fullref{lemmefinal}, so
$\smash{\mathrm{Sq}^{2^i}\mathrm{Sq}^{2^i}}w_{i,0}$ reduces to zero in $\mathrm{R}_0
\mathrm{F}_{-2}\cohom X_0$. Hence $\alpha _{i+1,0}\cup \alpha _{i+1,0}$
reduces to zero in $\mathrm{R}_0 \mathrm{F}_{-2}\cohom X_0$ for $i\geq\kappa$.

In other words, the element $\alpha _{i,0}\cup \alpha _{i,0}$ is
nilpotent for $i\geq\kappa$ and thus for some $t$,
$$\mathrm{Sq}_0^t (\alpha _{i,0}\cup \alpha _{i,0}) =\mathrm{Sq}_0^t\alpha
_{i,0}\cup \mathrm{Sq}_0^t\alpha _{i,0}=0.$$
This completes the proof that the cup square of $\alpha _{i,0}$ is trivial.

On the other hand, we know that 
$$\mathrm{Sq}_0^t\alpha _{i,0}\cup \mathrm{Sq}_0^t\alpha _{i,0}=\alpha
_{i+t,0}\cup \alpha _{i+t,0}=\mathrm{Sq}^{2^{i+t}}\alpha _{i+t,0}=\alpha
_{i+t+1,0}\neq 0.$$
This is a contradiction and completes the proof of the main theorem.

\appendix
\section{Trivial Bockstein actions and Lannes' functor}
The material of this section is well-known. It is already used in
\cite[Proposition 1.3, page 328]{Ku} and first proved by M\,Winstead \cite{W}.
We thank gratefully J\,Lannes who explained us the following proof.
Let $M$ be an unstable module. The notation $M^{\geq n}$ stands for the
submodule of $M$ of elements of degrees greater than $n$. We say that the
action of the Bockstein is \emph{trivial in degree greater than $n$\/} if
$\mathrm{Sq}^1 M^{\geq n}=0$.
\begin{proposition}
\label{appbock0}
Let $M$ be an unstable module. The action of the Bockstein in $M$ is trivial
in degree greater than $n$ if and only if  the action of the Bockstein in
$\mathrm{T} M$  is trivial in degrees greater than $n$.
\end{proposition}
Because $\redt M$ is a submodule of $\mathrm{T} M$ we have:
\begin{cor}
\label{appbock}
Let $M$ be an unstable module. If the action of the Bockstein in $M$ is
trivial in degree greater than $n$, then the action of the Bockstein in
$\redt M$  is also trivial in degrees greater than $n$.
\end{cor}
Before proving \fullref{appbock0}, we recall  the definition of the \emph{double\/}  $\Phi M$ of an
unstable module $M$ \cite[page 27]{LZ,Sc1}. The module $\Phi M$ is the unique unstable module $\Phi
M$ such that

\begin{enumerate}
\item[(i)] the module $\Phi M$ is zero in odd degrees,
\item[(ii)] for any $\ell$, $(\Phi M)^{2\ell}$ is $M^\ell$,
\item[(iii)] the natural map $\Phi\co \Phi M\lra M$ which maps $m$ to $\Phi m
=\mathrm{Sq}_0 m$ is linear with respect to the Steenrod algebra.
\end{enumerate}
$$\mathrm{Sq}^{2\ell}\Phi m =\Phi \mathrm{Sq}^\ell m.\leqno{\hbox{In other words}}$$
It is evident from the definition that the action of the Bockstein  is
trivial on $\Phi M$. Conversely, we have

\begin{lemme}
\label{lemdec}
Let $M$ be an unstable module such that the action of the Bockstein is
trivial in each degree. Denote by $M^{\text{odd}}$ and $M^{\text{even}}$ the odd and
even degree parts of $M$ as graded vector spaces. Then $M$ splits as a
module over the Steenrod algebra as
$$ M= M^{\text{odd}} \oplus M^{\text{even}}.$$
\end{lemme}
\begin{proof} This lemma is the consequence of the
following facts:

\begin{enumerate}
\item[(i)] the Steenrod algebra is generated as an algebra by the squares
$\mathrm{Sq}^{i}$,
\item[(ii)] we have for any odd square the Adem relation
$$\mathrm{Sq}^{2n+1}=\mathrm{Sq}^1 \mathrm{Sq}^{2n}. $$
\end{enumerate}

When the action of the Bockstein is trivial, it follows that
$M^{\text{odd}}$ and $M^{\text{even}}$ are unstable submodules and that the
vector space decomposition $M= M^{\text{odd}} \oplus M^{\text{even}}$ is in fact
a Steenrod algebra module decomposition. \fin
\begin{lemme}
\label{lemmedec2}
Let $M$ be a module such that $M$ is zero in odd degrees. Then $M$ is of the
form $\Phi M_1$ for a unique unstable module $M_1$. Let $M$ be an unstable
module such that $M$ is zero in even degrees. Then $M$ is of the form
$M=\Sigma \Phi M_2$ for a unique module $M_2$.
\end{lemme}
\begin{proof} Let us prove the first assertion. It
follows from the definitions that
$M_1$ has to be defined by $M_1^\ell = M^{2\ell}$. Furthermore, we also have
no choice for the Steenrod algebra structure on $M_1$. It remains only to
show that this actually defines an action of the Steenrod algebra, which
amounts to the definition of $\Phi$.

To prove the second assertion, we remark that for any module $M$
concentrated in odd degrees, the operator $\mathrm{Sq}_0$ is trivial. But
The triviality of this operator is exactly the obstruction for algebraically
desuspending an unstable module. So $M$ is of the form $M=\Sigma M'$ for a
unique $M'$. Now $M'$ is concentrated in even degree, and by the first part,
we have that $M'=\Phi M_2$ for a unique $M_2$. So, we have
$$M=\Sigma M'= \Sigma \Phi M_2.\proved$$
\fin

We return to the proof of  \fullref{appbock0}.

\begin{proof}[Proof of  \fullref{appbock0}]
Let $M$ be an unstable module having trivial action of the Bockstein in
degrees greater than $n$.

We have a short exact sequence of unstable modules
$$M^{\geq n} \lra M \lra M/M^{\geq n}.$$
By exactness of the $\mathrm{T}$ functor, we get an exact sequence
\begin{eqnarray}\label{one}
\mathrm{T}M^{\geq n} \lra \mathrm{T}M \lra \mathrm{T}(M/M^{\geq n}).
\end{eqnarray}
Lannes' $\mathrm{T}$ functor admits a natural splitting as
$$T\cong \redt \oplus  \mathrm{Id}$$
hence the exact sequence \eqref{one} splits into two short exact sequences
$$ M^{\geq n}\lra  M \lra   M   /M^{\geq n} \quad\text{and}\quad\redt M^{\geq
n}\lra \redt  M \lra   \redt M /\redt M^{\geq n}\cong \redt (M /M^{\geq n}).$$
Now $M/M^{\geq n}$ is a bounded module, so $\redt (M /M^{\geq n})=0$. On the
other hand, $M^{\geq n}$ has trivial action of Bocksteins and so, by \fullref{lemdec},
$$M^{\geq n}={(M^{\geq n})}^{\text{even}}\oplus {(M^{\geq n})}^{\text{odd}}.$$
Now, \fullref{lemmedec2} ensures that
\begin{gather*}M^{\geq n}=\Phi M_1 \oplus\Sigma \Phi M_2\\
\tag*{\hbox{and so}}
\redt M^{\geq n}=\redt (\Phi M _1\oplus \Sigma \Phi M _2) = (\Phi \redt M
_1 \oplus \Sigma \Phi \redt M _2)\end{gather*}
because the functor $\redt$ commutes to suspensions and to $\Phi$.

It follows
that $\redt M^{\geq n}$ has trivial action of Bocksteins in each degrees.
Finally, $\lt M= M\oplus \redt M$ has trivial action of Bocksteins in
degrees greater than $n$.

The converse is a consequence of the  aforementioned splitting of the
$\lt$ functor.
\fin

\bibliographystyle{gtart}
\bibliography{link}

\end{document}